# An insertion algorithm over staircase tableaux compatible with the ASEP's matrix ansatz.

Patxi Laborde-Zubieta[*][1]

[1]*LaBRI, University of Bordeaux*

**Abstract.** Based on the matrix ansatz of Derrida, Evans, Hakim and Pasquier, we prensent a new way of computing the stationary probability of a state of the asymmetric simple exclusion process (ASEP). Through an insertion algorithm over staircase tableaux, we give a combinatorial proof to the current interpretation of the ASEP by these tableaux of Corteel and Williams. The insertion algorithm induces a recursive structure which implies nice factorised formulas for the generating polynomials of staircase tableaux, as well as a bijection with some coloured inversion tables. In addition, we adapt the insertion algorithm to the case of type *B* symmetric tableaux and we define a new matrix ansatz compatible with it.

## 1 Introduction

The asymmetric simple exclusion process (ASEP) is a reference model of non-equilibrium statistical mechanics, it depends on six parameters $\alpha, \beta, \gamma, \delta, q, u \in [0,1]$. It describes, through a Markov chain, a system of particles hoping left and right on finite one dimensional lattice. Although it is simple, it exhibits non trivial macroscopic properties that we find generically in more realistic models. The seminal work of Corteel and Williams [6, 7, 8] initiated a fruitful combinatorial interpretation of the unique probability distribution of the ASEP with tableaux. Their interpretation is based on the matrix ansatz of Derrida, Evans, Hakim and Pasquier [10], which is a mathematical tool allowing the computation of the unique stationary distribution of the ASEP.

In the case $\gamma = \delta = 0$, the theory of the *cellular ansatz* of Viennot [13] shows that the tableaux involved in the combinatorial interpretation of the ASEP encode the computations of the stationary probabilities. While his approach gives a better understanding of the link between the matrix ansatz and the tableaux in the case $\gamma = \delta = 0$, it is difficult to apply it directly to the general case.

In this work, we follow the same idea as Viennot. Based on the matrix ansatz, we describe a new recursive way to compute the stationary probability of a state. Then, through an insertion algorithm over staircase tableaux, we show that generating polynomials of these tableaux follow the same recurrences. Up to an alternative proof of the existence of a solution to the generalised matrix ansatz (Remark 3), we get a combinatorial proof of the interpretation of the ASEP with staircase tableaux. The current proof of [8] is long and computational. As a consequence of the insertion algorithm, we obtain easily nice factorised formulas for the generating polynomial of staircase tableaux of fixed size, as well as a bijection with coloured inversion tables in which, contrary to the bijection of Corteel and Dasse-Hartaut [3], the parameter $q$ as an interpretation. This might be a way to achieve an interpretation of the general ASEP with coloured permutations.

---

[*]plaborde@labri.fr.

Type $B$ staircase tableaux [3] are symmetric staircase tableaux with respect to an involution motivated by the symmetries of the ASEP. In the last section, we adapt our insertion algorithm to the type $B$. Thus, we generalise factorised formulas of the generating function obtained by Corteel and Dasse-Hartaut [3]. Moreover, we define a new matrix ansatz compatible with the recurrences induced by the insertion algorithm, thus generalising the matrix ansatz of Corteel, Josuat-Vergès and Williams [4]. As in the case of staircase tableaux, we still miss the existence of a solution for this new matrix ansatz.

It should be noted that this work can be generalised to the two-species ASEP and rhombic staircase tableaux [5]. As of now, no preprint is available.

## 2 A new way of studying the matrix ansatz

Let $n$ be a positive integer. The states of the ASEP are the words of size $n$ with $\{\bullet, \circ\}$ as alphabet. The symbol $\circ$ represents an empty site and $\bullet$ a site occupied by a particle. We also use the terminology of a white and black particle. The possible transitions from a state are illustrated in Figure 1. More precisely, the transition probability $p(s_1, s_2)$ between two states $s_1$ and $s_2$ is defined as follows ($X$ and $Y$ represent words in $\circ$ and $\bullet$): if $s_1 = X\bullet\circ Y$ and $s_2 = X\circ\bullet Y$, then $p(s_1, s_2) = \frac{u}{n+1}$ and $p(s_2, s_1) = \frac{q}{n+1}$; if $s_1 = X\bullet$ and $s_2 = X\circ$, then $p(s_1, s_2) = \frac{\beta}{n+1}$ and $p(s_2, s_1) = \frac{\delta}{n+1}$; if $s_1 = \circ Y$ and $s_2 = \bullet Y$, then $p(s_1, s_2) = \frac{\alpha}{n+1}$ and $p(s_2, s_1) = \frac{\gamma}{n+1}$; in any other case if $s_1 \neq s_2$, then $p(s_1, s_2) = 0$; finally $p(s_1, s_1) = 1 - \sum_{s_2 \neq s_1} p(s_1, s_2)$. For general parameters, the ASEP has a unique stationary distribution. Let $n$ be a positive integer, we denote $\mu_n(m_n \cdots m_1)$ the stationary probability of a state $m_n \cdots m_1 \in \{\bullet, \circ\}^n$. Let us consider the notations $\widetilde{\circ} = \bullet$ and $\widetilde{\bullet} = \circ$. Moreover, the symbol $Q|_{x \leftrightarrow y}$ means that the parameters $x$ and $y$ are interchanged in $Q$. The ASEP has the following symmetries which induce identities satisfied by $\mu_n$.

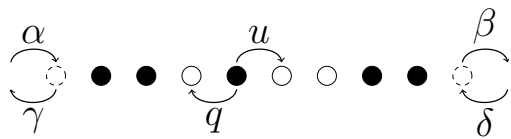

**Figure 1:** The possible transitions of the ASEP.

(C) The *C-symmetry* (charge-conjugation) consists in replacing each $p(s, s')$ with $p(s', s)$ in the previous definition. By doing so, we obtain the same model up to the exchange of the parameters $\alpha$ and $\gamma$ on the one hand, $\beta$ and $\delta$ on the other hand and $q$ and $u$ lastly. As a consequence $\mu_n(m_n \cdots m_1) = \mu_n(\widetilde{m}_n \cdots \widetilde{m}_1)|_{\alpha \leftrightarrow \gamma, \beta \leftrightarrow \delta, q \leftrightarrow u}$.

(P) The *P-symmetry* (parity symmetry) consists in replacing each $p(s,s')$ with $p(\bar{s},\bar{s'})$, where $\overline{m_n \cdots m_1} = m_1 \cdots m_n$. By doing so, we obtain the same model up to the exchange of the parameters $\alpha$ and $\delta$ on the one hand, $\beta$ and $\gamma$ on the other hand and $q$ and $u$ lastly. It follows that $\mu_n(m_n \cdots m_1) = \mu_n(m_1 \cdots m_n)|_{\alpha \leftrightarrow \delta, \beta \leftrightarrow \gamma, q \leftrightarrow u}$.

(CP) The *CP-symmetry* is the combination of the two previous symmetries. We obtain the same model up to the exchange of the parameters $\alpha$ and $\beta$ on the one hand and $\gamma$ and $\delta$ on the other hand. It implies that $\mu_n(m_n \cdots m_1) = \mu_n(\widetilde{m_1} \cdots \widetilde{m_n})|_{\alpha \leftrightarrow \beta, \gamma \leftrightarrow \delta}$.

The matrix ansatz of Derrida *et al.* introduced in [10] is a system of equations which allows to compute the stationary distribution of the ASEP. If we find two matrices, $D$ and $E$, and two vectors, $V$ and $W$, satisfying the matrix ansatz, then the stationary probability of any state is expressed as the product $WXV$, where $X$ is a product in $D$ and $E$ depending on the state. Corteel and Williams generalise in [8] the matrix ansatz of Derrida *et al.* by adding a degree of freedom. As a result, they are able to construct a solution without any parameter restrictions, while there exists no matrix solution to the initial matrix ansatz when $\alpha\beta = q^i\gamma\delta$ ([11, Section IV]). The equations of the matrix ansatz give also a recursive way to compute the stationary probability of a state without any matrix representation. As in [2] but in a different way, we study the matrix ansatz as a linear system with an infinite number of equations. Let $\lambda = (\lambda_n)_{n \geq 0} \in \mathbb{Q}(\alpha, \beta, \gamma, \delta, q, u)^{\mathbb{N}}$ and $g^\lambda : \{\bullet, \circ\}^* \to \mathbb{Q}(\alpha, \beta, \gamma, \delta, q, u)$ a function solution of

$$\begin{cases} ug^\lambda(X\bullet\circ Y) = qg^\lambda(X\circ\bullet Y) + \lambda_{|X|+|Y|+2}\left(g^\lambda(X\bullet Y) + g^\lambda(X\circ Y)\right) & (I) \\ \alpha g^\lambda(\circ X) = \gamma g^\lambda(\bullet X) + \lambda_{|X|+1}g^\lambda(X) & (II) \\ \beta g^\lambda(X\bullet) = \delta g^\lambda(X\circ) + \lambda_{|X|+1}g^\lambda(X) & (III) \\ g^\lambda(\varepsilon) = \lambda_0, \end{cases} \quad (2.1)$$

where $|X|$ is the size of $X$ and $\varepsilon$ the empty word. The original case corresponds to $\lambda = 1$. The stationary probability of a state is obtained as follows.

**Theorem 1** ([8, Theorem 5.2]). *Let* $m = m_n \cdots m_1 \in \{\bullet, \circ\}^n$. *Assume that for a sequence* $\lambda = (\lambda_n)_{n \geq 0} \in \mathbb{Q}(\alpha, \beta, \gamma, \delta, q, u)^{\mathbb{N}}$, *there exists a solution* $g^\lambda$ *to the system of equations* (2.1), *then* $\mu_n(m) = \frac{g^\lambda(m)}{Z_{n,g^\lambda}(\alpha, \beta, \gamma, \delta; q, u)}$, *where* $Z_{n,g^\lambda}(\alpha, \beta, \gamma, \delta; q, u) = \sum_{m' \in \{\bullet, \circ\}^n} g^\lambda(m')$.

The generating function $Z_{n,g^\lambda}$ is called the *partition function*. From now on, the parameters $\alpha, \beta, \gamma, \delta, q$ and $u$ are viewed as indeterminates.

Let $m$ be a state of size $n+1$. Let $\circ(m)$ and $\bullet(m)$ denote respectively the number of white and black particles of $m$. We can compute $g^\lambda(m)$ by moving the leftmost particle of $m$ all the way to the right and then all the way to the left. In order to lighten the notations, we write $m'$ instead of $g^\lambda(m')$. Let us suppose that $m = \bullet X$, we get

$$\begin{aligned}
m &= \bullet X \\
&= \left(\frac{q}{u}\right)^{|\circ(X)|} X\bullet + \sum_{X=U\circ V} \frac{\lambda_{n+1}}{u}\left(\frac{q}{u}\right)^{|\circ(U)|}(U\bullet V + U\circ V) \\
&= \left(\frac{q}{u}\right)^{|\circ(X)|} \frac{\delta}{\beta} X\circ + \left(\frac{q}{u}\right)^{|\circ(X)|} \frac{\lambda_{n+1}}{\beta} X + \sum_{X=U\circ V} \frac{\lambda_{n+1}}{u}\left(\frac{q}{u}\right)^{|\circ(U)|}(U\bullet V + U\circ V) \\
&= \left(\frac{q}{u}\right)^{n} \frac{\delta}{\beta} \circ X + \frac{\delta}{\beta} \sum_{X=U\bullet V} \frac{\lambda_{n+1}}{u}\left(\frac{q}{u}\right)^{|\circ(U)|+|V|}(U\bullet V + U\circ V)+ \\
&\quad \left(\frac{q}{u}\right)^{|\circ(X)|} \frac{\lambda_{n+1}}{\beta} X + \sum_{X=U\circ V} \frac{\lambda_{n+1}}{u}\left(\frac{q}{u}\right)^{|\circ(U)|}(U\bullet V + U\circ V) \\
&= \left(\frac{q}{u}\right)^{n} \frac{\gamma\delta}{\alpha\beta} m + \left(\frac{q}{u}\right)^{n} \frac{\delta}{\alpha\beta}\lambda_{n+1} X + \frac{\delta}{\beta}\sum_{X=U\bullet V} \frac{\lambda_{n+1}}{u}\left(\frac{q}{u}\right)^{|\circ(U)|+|V|}(U\bullet V + U\circ V)+ \\
&\quad \left(\frac{q}{u}\right)^{|\circ(X)|} \frac{\lambda_{n+1}}{\beta} X + \sum_{X=U\circ V} \frac{\lambda_{n+1}}{u}\left(\frac{q}{u}\right)^{|\circ(U)|}(U\bullet V + U\circ V).
\end{aligned}$$

Therefore, we obtain the first (global) recurrence (equation 2.2) of Theorem 2. Similarly, if $m = \circ X$ we obtain the second (global) recurrence (equation 2.3). We notice that $g^\lambda(m)$ is a linear combination of weights of smaller words. If the sequence $(\lambda_n)_{n\geqslant 0}$ is never equal to zero, then $g^\lambda$ is entirely determined by $g^\lambda(\varepsilon) = \lambda_0$ and the equations (2.2) and (2.3) of Theorem 2. Hence, the system (2.1) has at most one solution. It remains to prove that a function defined by these two equations is indeed a solution to the matrix ansatz.

**Theorem 2.** *Let $\lambda = (\lambda_n)_{n\geqslant 0} \in \mathbb{Q}(\alpha,\beta,\gamma,\delta,q,u)^{\mathbb{N}}$ be a sequence. The function $h^\lambda : \{\bullet,\circ\}^* \to \mathbb{Q}(\alpha,\beta,\gamma,\delta,q,u)$ defined by the equations $h^\lambda(\varepsilon) = \lambda_0$,*

$$h^\lambda(\bullet X) = \frac{\lambda_{n+1}}{u^n\alpha\beta - q^n\gamma\delta}\left[\alpha\delta \sum_{X=U\bullet V} u^{n-1-|\circ(U)|-|V|}q^{|\circ(U)|+|V|}(h^\lambda(U\bullet V)+h^\lambda(U\circ V))+\right.$$
$$\left. q^n\delta X + u^{n-|\circ(X)|}q^{|\circ(X)|}\alpha X + \alpha\beta\sum_{X=U\circ V} u^{n-1-|\circ(U)|}q^{|\circ(U)|}(h^\lambda(U\bullet V)+h^\lambda(U\circ V))\right] \quad (2.2)$$

*and*
$$h^\lambda(\circ X) = \frac{\lambda_{n+1}}{u^n\alpha\beta - q^n\gamma\delta}\left[\delta\gamma \sum_{X=U\bullet V} u^{n-1-|\circ(U)|-|V|}q^{|\circ(U)|+|V|}(h^\lambda(U\bullet V)+h^\lambda(U\circ V))+\right.$$
$$\left. u^{n-|\circ(X)|}q^{|\circ(X)|}\gamma X + \gamma\beta\sum_{X=U\circ V} u^{n-1-|\circ(U)|}q^{|\circ(U)|}(h^\lambda(U\bullet V)+h^\lambda(U\circ V)) + u^n\beta X\right] \quad (2.3)$$

*is a solution of the system 2.1.*

*Sketched proof.* If the theorem is true for a specific sequence which is never equal to zero, then it is true for all sequences. Moreover, it is sufficient to prove that the system has a solution, since this solution would be necessarily equal to $h^\lambda$. Yet, for $\lambda_0 = 1$ and $\lambda_n = u^{n-1}\alpha\beta - q^{n-1}\gamma\delta$, Corteel and Williams [8] exhibits a solution. □

**Remark 3.** *The purpose of our approach is to provide a combinatorial proof to the interpretation of the stationary distribution of the ASEP by staircase tableaux. The current proof of Corteel and Williams is long and computational. Through an insertion algorithm over staircase tableaux, we prove in Theorem 9 that $h^\lambda(m)$ have a combinatorial interpretation as a generating polynomial of staircase tableaux. Hence Theorems 2 and 9 provide the desired combinatorial interpretation. However, for the moment the proof of Theorem 2 uses the result that we want to prove.*

*A proof using only the recurrence relations seems within reach. The equation (I) of system 2.1 is a direct consequence. The equation (III) requires more work. As for the last one, we are not able to prove it for the moment. An alternative demonstration could be to study the system 2.1 as a linear system and to prove the existence of a solution by computing the row echelon form of the associated matrix.*

*As it is, Theorem 2 provides a new way to show that a combinatorial class, different from staircase tableaux, gives an interpretation of the stationary distribution of the ASEP.*

Let $w$ be the function $h^\lambda$ with $\lambda_n = u^{n-1}\alpha\beta - q^{n-1}\gamma\delta$ for all $n \geqslant 1$ and $\lambda_0 = 1$. The value $w(m)$ is called the *weight of $m$*. The advantage of this choice of $\lambda$ is that $w$ has value in $\mathbb{Q}[\alpha, \beta, \gamma, \delta, q, u]$, more precisely, if $m$ is of size $n$, the $w(m)$ is a positive sum of $2^n n!$ monomials of degree $\frac{n(n+1)}{2}$. In order to see that, we simply have to notice that $w(m)$ is a positive linear combination of the weights of $2n$ words of size $n-1$ whose coefficients are monomials of degree $n$. In particular, we get $Z_{n,w}(1,1,1,1;1,1) = 4^n n!$. Since the degree of the monomials of $w(m)$ is constant equal to $\frac{n(n+1)}{2}$, in the rest of the article, we assume without loss of generality that $u = 1$.

Let $m \in \{\bullet, \circ\}^n$, the first column of the Table 1 describes the words $m'$ of size $n+1$ such that $w(m)$ appears in the recurrence identity, (2.2) or (2.3), satisfied by $w(m')$. The second column gives the associated coefficient. The arrows, under the words $m'$, illustrate the course of the leftmost particle resulting in $m$.

If we specialise $q$ to 1, the Table 1 becomes simpler, so that $Z_{n+1,w}(\alpha, \beta, \gamma, \delta; 1, 1)$ is equal to

$$(\alpha + \beta + \gamma + \delta + n(\alpha + \delta)(\beta + \gamma))Z_{n,w}(\alpha, \beta, \gamma, \delta; 1, 1).$$

| State | Coefficient |
|---|---|
| $\bullet\, m \rightarrow$ | $q^{|\circ(m)|}\alpha$ |
| $\bullet\, m \rightleftarrows$ | $q^n \delta$ |
| $\circ\, m \hookrightarrow$ | $q^{|\circ(m)|}\gamma$ |
| $\circ\, m \leftarrow$ | $\beta$ |
| For every decomposition $m = U\bullet V$ or $m = U\circ V$ | |
| $\bullet\, U \bullet V \rightleftarrows$ | $q^{|\circ(U)|+|V|}\alpha\delta$ |
| $\bullet\, U \circ V \rightarrow$ | $q^{|\circ(U)|}\alpha\beta$ |
| $\circ\, U \bullet V \rightleftarrows$ | $q^{|\circ(U)|+|V|}\gamma\delta$ |
| $\circ\, U \circ V \hookrightarrow$ | $q^{|\circ(U)|}\gamma\beta$ |

Table 1: Preimages of $m \in \{\bullet, \circ\}^n$.

As a consequence, we get $Z_{n,w}(\alpha, \beta, \gamma, \delta; 1, 1) = \prod_{i=0}^{n-1}(\alpha + \beta + \gamma + \delta + i(\alpha + \delta)(\beta + \gamma))$. This identity has several combinatorial proofs based on staircase tableaux, or equivalent tableaux, [9, 1, 3]. In the original case where $\lambda$ is the constant sequence equal to 1, the equivalent formula was proved in [12, 2]. Our proof is the most elementary one.

## 3　An insertion algorithm over staircase tableaux

Staircase tableaux were introduced by Corteel and Williams in [8].

**Definition 4** (Staircase tableau). *A staircase tableau of size n is a Ferrers diagram of shape $(n, n-1, \ldots, 1)$ such that cells are either empty or labelled with a symbol $\{\alpha, \beta, \gamma, \delta\}$ while respecting the rules: each cell of the south-est border is labelled; the cells above $\alpha/\gamma$, i.e. $\alpha$ or $\gamma$, are empty; the cells to the left of $\beta/\delta$ are empty.*

An example is given in the Figure 2. We denote by $\mathcal{TE}$ the set of staircase tableaux and $\mathcal{TE}_n$ those of size $n$. The weight $w(T)$ of a staircase tableau $T$ is equal to the product of all its symbols after labelling the empty cells as follows. An empty cell to the left of $\delta$ gets a $q$. It is also the case of an empty cell to the left of a $\alpha/\gamma$ and above $\beta/\gamma$. The remaining empty cells are labelled with $u$. The generating function of staircase tableaux of size $n$ with respect to this weight is denoted by $Z_n(\alpha, \beta, \gamma, \delta; q, u)$. The weight of a staircase tableau of size $n$ is of total degree $n(n+1)/2$, hence, we assume without loss of generality that $u = 1$. Let $T$ be a staircase tableau of size $n$. We label from 1 to $n$ the rows (resp. columns) bottom to top (resp. left to right). A cell is therefore uniquely identified by the pair $(i, j)$ where $i$ (resp. $j$) is the line (resp. column) number $i$ (resp. $j$). We denote $L_i$, $C_j$, $c_{i,j}$, $c_i$ respectively the line $L$ number $i$, the column number $j$, the cell of coordinates $(i, j)$ and the cell $c_{i,i}$. A row $L_i$ is of *type* $\alpha/\gamma$ if $c_i$ is labelled with $\alpha/\gamma$. Let $m = m_n \cdots m_1$ be the word such that $m_i = \bullet$ if $c_i = \alpha/\delta$ and $m_i = \circ$ if $c_i = \beta/\gamma$. The state $m$ is the *type* of $T$ and we write $type(T) = m$.

**Figure 2:** A staircase tableau of size 8, weight $u^8 q^{14} \alpha^3 \beta^4 \gamma^3 \delta^4$ and type $\bullet\circ\bullet\bullet\circ\bullet\bullet\circ$.

In order to introduce the insertion algorithm, we define some elementary operations over staircase tableaux without the filling of $q$ and $u$. Let $i_1 < i_2$ such that all the cells (resp. all the cells but $c_{i_2}$) of $L_{i_2}$ are empty. The operation $L_{i_1} \to L_{i_2}$ (resp. $L_{i_1}^* \to L_{i_2}$) moves the label of $c_{i_1,k}$ to $c_{i_2,k}$ for $1 \leqslant k \leqslant i_1$ (resp. $1 \leqslant k < i_1$). Finally, $T + L_{n+1}(x)$ adds an empty row of size $n + 1$ above $T$ and labels its rightmost cell with $x$.

**Definition 5.** *Let $T$ be a staircase tableau of size n. We define $4(n+1)$ choices of insertions. An example of the three cases is given in Figure 3*

- *First case. We insert a row with $\beta/\delta$: $T + L_{n+1}(y)$ with $y \in \{\beta, \delta\}$. We denote by $insertion(T, y)$ the resulting tableau.*

- *Second case. We start by adding a row with $\alpha/\gamma$: $T + L_{n+1}(x)$ with $x \in \{\alpha, \gamma\}$. Afterwards, we shift upwards the labels, except the rightmost one, of the rows of type $\alpha/\gamma$: let $i_1 < \cdots < i_k = n+1$ be the indices of the rows of type $\alpha/\gamma$ of $T + L_{n+1}(x)$. Then, we perform the operations $L_{i_{k-1}}^* \to L_{i_k}, \cdots, L_{i_1}^* \to L_{i_2}$. We denote by $insertion(T, x)$ the resulting tableau.*

- *Third case. We start by adding a row with $\alpha/\gamma$: $T + L_{n+1}(x)$ with $x \in \{\alpha, \gamma\}$. Afterwards, we choose $i \in [\![1, n]\!]$ and we shift upwards the labels, except the rightmost one, of the rows of type $\alpha/\gamma$ strictly above $L_i$. Let $i < i_l < \cdots < i_k = n+1$ be the indices of the rows of type $\alpha/\gamma$ of $T + L_{n+1}(x)$ higher than $L_i$. Then, we perform the operations $L^*_{i_{k-1}} \to L_{i_k}, \cdots, L^*_{i_l} \to L_{i_{l+1}}$ Finally, we move the labels of $L_i$ to $L_{i_l}$ and we label with $y \in \{\beta, \gamma\}$ the cell $c_i$: $L_i \to L_{i_l} T$ and $c_i = y$. We denote by $insertion(T, (x, y, i))$ the resulting tableau.*

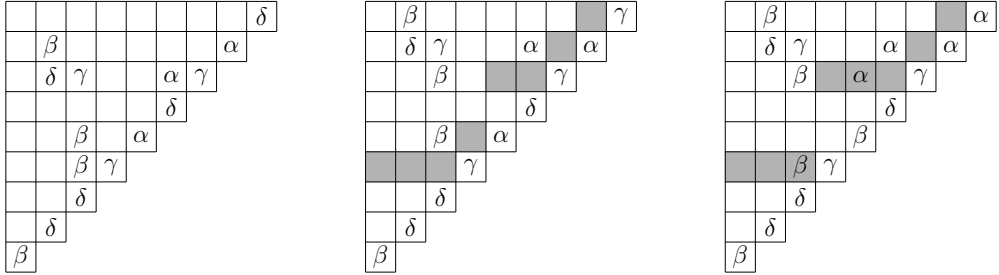

**Figure 3:** From left to right: $insertion(T, \delta)$, $insertion(T, \gamma)$ and $insertion(T, (\alpha, \beta, 5))$, where $T$ is the staircase tableau of Figure 2. The special cells are in grey.

The insertion algorithm endows staircase tableaux with a recursive structure.

**Theorem 6.** *The function insertion is a bijection between the set of pairs $(T, e)$, where $T$ is a staircase tableau of size $n$ and $e$ is an element among $\{\alpha, \beta, \gamma, \delta\} \cup \{(x, y, i) \mid x \in \{\alpha, \gamma\}, y \in \{\beta, \delta\}$ and $i \in [\![1, n]\!]\}\}$, and the set of staircase tableaux of size $n+1$.*

*Sketched proof.* Let $y \in \{\beta, \gamma\}$. Staircase tableaux of size $n+1$ such that $c_{n+1}$ is labelled with $y$ are trivially in bijection with pairs $(T, y)$. Let us now consider staircase tableaux of size $n+1$ such that $c_{n+1}$ is labelled with $x \in \{\alpha, \gamma\}$. We denote by $i_1 < \cdots < i_k = n+1$ the indices of the rows of type $\alpha/\gamma$. We call *special cells* the cells $\{c_{i_p, j} \mid p \in [\![2, k]\!], i_{p-1} \leq j < i_p\} \cup \{c_{i_1, j} \mid j \in [\![1, i_1 - 1]\!]\}$. They are in grey in the Figure 3. The staircase tableaux of size $n+1$ with all the special cells empty are obtained by an insertion of the form $insertion(T, x)$. The staircase tableaux of size $n+1$ such that the north-eastern most non empty special cell belongs to $C_i$ are obtained by an insertion of the form $insertion(T, (x, y, i))$ where $y$ is the label of $c_i$. □

A recursive structure was also defined in [1] but $q$ was difficult to track.

**Proposition 7.** *Let $T$ be a staircase tableau of size $n$, $w = w(T)$ and $m = type(T)$. For each of the following cases, $m'$ refers to the type of the staircase tableau obtained after the insertion.*

- $w(insertion(T, \alpha)) = q^{|\circ(m)|} \alpha w$ and $m' = \bullet m$,
- $w(insertion(T, \beta)) = \beta w$ and $m' = \circ m$,
- $w(insertion(T, \gamma)) = q^{|\circ(m)|} \gamma w$ and $m' = \circ m$,
- $w(insertion(T, \delta)) = q^n \delta w$ and $m' = \bullet m$.

*For every decomposition* $m = U \overset{i}{\bullet} V$ *or* $m = U \overset{i}{\circ} V$ *we have,*

- $w(insertion(T, (\alpha, \beta, i))) = q^{|\circ(U)|}\alpha\beta w$ *and* $m' = \bullet U \circ V$,
- $w(insertion(T, (\gamma, \beta, i))) = q^{|\circ(U)|}\gamma\beta w$ *and* $m' = \circ U \circ V$,
- $w(insertion(T, (\alpha, \delta, i))) = q^{|\circ(U)|+|V|}\alpha\delta w$, *and* $m' = \bullet U \bullet V$,
- $w(insertion(T, (\gamma, \delta, i))) = q^{|\circ(U)|+|V|}\gamma\delta w$, *and* $m' = \circ U \bullet V$.

*Sketched proof.* For $w(insertion(T, y))$ with $y \in \{\beta, \delta\}$ the proposition is immediate. Else, if we fill $T$ with $q$ and $u$, we observe that the operations $L^*_{i_{j-1}} \to L_{i_j}$ and $L_{i_{j-1}} \to L_{i_j}$ of the insertion algorithm moves also the symbols $q$ and $u$. Hence, we have to identify the new $n+1$ symbols. In the following explanation, $i$ equals 0 for $insertion(T, x)$ with $x \in \{\alpha, \gamma\}$. The new symbols are in the special cells in the columns $C_j$ for $j > i$ and the cells $c_{j,i}$ for $j \leqslant i$. The symbol of a special cell $c_{j,i}$ is determined by the type, $\circ$ or $\bullet$, of $c_i$. The symbols of the cells $c_{j,i}$ for $j \leqslant i$ depend only on the symbol $y$ of the insertion. □

We get the factorised formula of $Z_n(\alpha, \beta, \gamma, \delta; 1, 1)$ and two identities of [3].

**Corollary 8.** *We have the identities*

$$Z_n(\alpha, \beta, \gamma, \delta; 1, 1) = \prod_{i=0}^{n-1}(\alpha + \beta + \gamma + \delta + i(\alpha + \delta)(\beta + \gamma))$$

$$Z_n(\alpha, 0, 0, \delta; q, u) = \prod_{i=0}^{n-1}(u^i\alpha + \alpha\delta(u^{i-1} + u^{i-2}q + \cdots + q^{i-1}) + q^i\delta).$$

$$Z_n(0, \beta, \gamma, 0; q, u) = \prod_{i=0}^{n-1}(u^i\beta + \beta\gamma(u^{i-1} + u^{i-2}q + \cdots + q^{i-1}) + q^i\gamma).$$

The weight of a state of the ASEP can now be interpreted with staircase tableaux.

**Theorem 9.** *Let* $n \geqslant 0$ *and* $m \in \{\bullet, \circ\}^n$, *then* $w(m) = \sum_{\substack{T \in \mathcal{TE}_n \\ type(T) = m}} w(T)$.

*Sketched proof.* It follows from the fact that proposition 7 is compatible with Table 1. □

The insertion algorithm over staircase tableaux gives a simple bijection between staircase tableaux and coloured inversion tables. Indeed, for a staircase tableau $T$ of size $n$, there is a unique sequence $T_1, \ldots, T_n$ of staircase tableaux such that $T_k$ is obtained from $T_{k-1}$ with an insertion. We denote by $e_1, \ldots, e_n$ the sequence such that $e_1$ is equal to the unique symbol of $T_1$ and $T_k = insertion(T_{k-1}, e_k)$. Let $I(T) = (i_1, \cdots, i_n)$ be the coloured inversion table such that $i_k = 0_x$ if $e_k = x \in \{\alpha, \gamma\}$, $i_k = k_y$ if $e_k = y \in \{\beta, \delta\}$ and $i_k = j_{x,y}$ if $e_k = (x, y, j)$. We just defined a bijection between staircase tableaux and some coloured inversion tables. With some more work, we can read the symbols of the south-east border of $T$ directly on $I(T)$. Moreover, it is not hard to define a statistic over these coloured inversion tables corresponding to the parameter $q$. Corteel and Dasse-Hartaut [3] defined a bijection between staircase tableaux and some coloured inversion tables, but the parameter $q$ could not be read directly on the inversion tables. In the case $\gamma = \delta = 0$, the statistic corresponding to $q$ we obtain is $\sum_{1 \leqslant k \leqslant n} |\{i_j, j < k \text{ and } i_j > i_k\}|$, while in the case $\alpha = \delta = 0$, it corresponds to $\sum_{1 \leqslant i \leqslant n} max(i - 1 - I[i], 0)$. Through this new bijection we can hope to achieve a bijection between staircase tableaux and some coloured permutation interpreting all the parameters as well as the symbols of the south-east border.

## 4 Staircase tableaux of type $B$

Type $B$ staircase tableaux were initially defined in [3]. We denote by $\psi$ the involution consisting of an axial symmetry with respect to the main diagonal and the exchange of $\alpha$ and $\delta$ on the one hand and $\beta$ and $\gamma$ on the other hand. It corresponds to the P-symmetry if we set $q = u = 1$.

**Definition 10.** *A* type $B$ staircase tableau *of size $n$ is a staircase tableau of size $2n$ invariant with respect to $\psi$.*

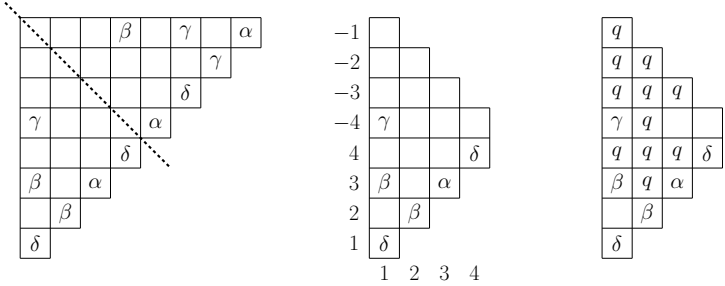

**Figure 4:** From left to right : a type $B$ staircase tableau of size 4, the corresponding half-tableau, the filling with $q$ of the half tableau.

All this section can be adapted to the tableaux corresponding to the CP-symmetry. A type $B$ staircase tableau is entirely defined by the bottom half. From now on, we work with bottom half tableaux. The type (resp. weight) of a of a type $B$ staircase tableau is the type (resp. weight) of its bottom half. We add an indeterminate $z$ which counts the number of $q$ in the main diagonal. The weight of the type $B$ staircase tableau of Figure 4 is $\alpha \beta^2 \gamma \delta^2 q^{11} u^3 z^3$ and its type is ●●○●. The generating polynomial of type $B$ staircase tableaux of size $n$ with respect to this weight is denoted $Z_n^B(\alpha, \beta, \gamma, \delta; q, u, z)$. If $z = 1$, the weight of a type $B$ staircase tableau of size $n$ is of degree $n(n+1)$. Hence, without loss of generality we can assume that $u = 1$.

As for staircase tableaux, in order to introduce the insertion algorithm over type $B$ staircase tableaux, we define elementary operations. Let $T$ be a type $B$ staircase tableau of size $n$. We number with $1, \ldots, n, -n, \ldots, -1$ (resp. $1, \ldots, n$) the rows (resp. columns) from bottom to top (resp. left to right). We change the operation $T + L_{n+1}(x)$ to $T +$ $L_{n+1}(x) + L_{-n-1}$ which consists in inserting two empty rows of size $n+1$ between $L_n$ and $L_{-n}$ and labelling the rightmost cell of $L_{n+1}$ with $x$.

**Definition 11.** *Let $T$ be a type B staircase tableau of size $n$, we define $4(2n+1)$ choices of insertions. The four cases are illustrated in Figure 5.*

- *First case. We start with $T + L_{n+1}(x) + L_{-n-1}$ with $x \in \{\alpha, \gamma\}$. Let $i_1 < \cdots < i_k = n+1$ be the indices of the rows of type $\alpha/\gamma$. We do $L^*_{i_{k-1}} \to L_{i_k}, \cdots, L^*_{i_1} \to L_{i_2}$. We denote insertion$(T, x)$ the tableau we obtain.*

- *Second case. We start with $T + L_{n+1}(x) + L_{-n-1}$ and $x \in \{\alpha, \gamma\}$. Let $i \in [\![1, n]\!]$ and $i < i_l < \cdots < i_k = n+1$ be the indices of the rows of type $\alpha/\gamma$ above $L_i$. Then, we perform $L^*_{i_{k-1}} \to L_{i_k}, \cdots, L^*_{i_l} \to L_{i_{l+1}}$. Finally, we do $L_i \to L_{i_l}$ and $c_i = y$ with $y \in \{\beta, \delta\}$. We denote insertion$(T, (x, y, i))$ the tableau we obtain.*

- *Third case. We start with $T + L_{n+1}(x) + L_{-n-1}$ and $x \in \{\beta, \delta\}$. Let $0 < i_1 < \cdots < i_{k-1}$ be the indices of the rows of type $\alpha/\gamma$ and $i_k = -n-1$. Then, we do $L^*_{i_{k-1}} \to L_{i_k}, \cdots, L^*_{i_1} \to L_{i_2}$. We denote insertion$(T, x)$ the tableau we obtain.*

- *Fourth case. We start with $T + L_{n+1}(x) + L_{-n-1}$ with $x \in \{\beta, \delta\}$. Let $i \in [\![1, n]\!]$ $i < i_l < \cdots < i_{k-1}$ be the indices of the rows of type $\alpha/\gamma$ above $L_i$ and $i_k = -n-1$. Then, we perform $L^*_{i_{k-1}} \to L_{i_k}, \cdots, L^*_{i_l} \to L_{i_{l+1}}$. Finally, we do $L_i \to L_{i_l}$ and $c_i = y$ with $y \in \{\beta, \delta\}$ We denote insertion$(T, (x, y, i))$ the tableau we obtained.*

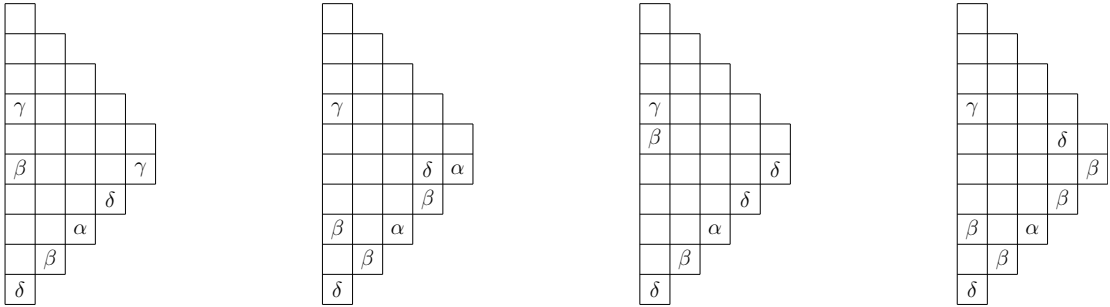

**Figure 5:** From left to right: *insertion*$(T, \gamma)$, *insertion*$(T, (\alpha, \beta, 4))$, *insertion*$(T, \delta)$ and *insertion*$(T, (\beta, \beta, 4))$, where $T$ is the type $B$ staircase tableau of Figure 4.

**Theorem 12.** *The function insertion is a bijection between the set of pairs $(T, e)$, where $T$ is a type B staircase tableau of size $n$ and $e$ an element among $\{\alpha, \beta, \gamma, \delta\} \cup \{(x, y, i) |\ x \in \{\alpha, \beta, \gamma, \delta\}, y \in \{\beta, \delta\}, i \in [\![1, n]\!]\}$, and the set of type B staircase tableaux of size $n+1$.*

We can describe the evolution of the weight depending on the type.

**Proposition 13.** *Let $T$ be a type B staircase tableau of size $n$, $w = w(T)$ and $m = type(T)$. For each of the cases, $m'$ refers to the type of the staircase tableau obtained after the insertion.*

- $w(insertion(T, \alpha)) = zq^{|\circ(m)|+n+1} \alpha w$ and $m' = \bullet m$,

- $w(insertion(T, \beta)) = zq^{|\circ(m)|+1}\, \beta w$ and $m' = \circ m$,
- $w(insertion(T, \gamma)) = q^{|\circ(m)|}\, \gamma w$ and $m' = \circ m$,
- $w(insertion(T, \delta)) = q^{|\circ(m)|+n}\, \delta w$ and $m' = \bullet m$.

*For any decomposition $m = U \overset{i}{\bullet} V$ or $m = U \overset{i}{\circ} V$, we have*

- $w(insertion(T, (\alpha, \beta, i))) = zq^{|\circ(U)|+n+1}\, \alpha\beta w$ and $m' = \bullet U \circ V$,
- $w(insertion(T, (\beta, \beta, i))) = zq^{|\circ(U)|+1}\, \beta\beta w$ and $m' = \circ U \circ V$,
- $w(insertion(T, (\gamma, \beta, i))) = q^{|\circ(U)|}\, \gamma\beta w$ and $m' = \circ U \circ V$,
- $w(insertion(T, (\delta, \beta, i))) = q^{|\circ(U)|+n}\, \delta\beta w$ and $m' = \bullet U \circ V$,
- $w(insertion(T, (\alpha, \delta, i))) = zq^{|\circ(U)|+|V|+n+1}\, \alpha\delta w$ and $m' = \bullet U \bullet V$,
- $w(insertion(T, (\beta, \delta, i))) = zq^{|\circ(U)|+|V|+1}\, \beta\delta w$ and $m' = \circ U \bullet V$,
- $w(insertion(T, (\gamma, \delta, i))) = q^{|\circ(U)|+|V|}\, \gamma\delta w$ and $m' = \circ U \bullet V$,
- $w(insertion(T, (\delta, \delta, i))) = q^{|\circ(U)|+|V|+n}\, \delta\delta w$ and $m' = \bullet U \bullet V$.

As a consequence, we generalise identities of Corteel and Dasse-Hartaut [3].

**Corollary 14.** *We have the identities*

$$Z_n^B(\alpha, \beta, \gamma, \delta; 1, 1, z) = (z(\alpha + \beta) + \gamma + \delta)^n \prod_{i=1}^{n-1}(1 + i(\beta + \delta))$$

$$Z_n^B(0, \beta, \gamma, ; q, u, z) = (zq\beta + u\gamma)^n \prod_{k=0}^{n-1}(q^k u^k + q^k(u^{k-1} + u^{k-2}q + \cdots + q^{k-1})\beta)$$

$$Z_n^B(\alpha, 0, 0, \delta; q, u, z) = (zq\alpha + u\delta)^n \prod_{k=0}^{n-1}(q^k u^k + u^k(u^{k-1} + u^{k-2}q + \cdots + q^{k-1})\delta).$$

In the case of staircase tableaux, the study of the matrix ansatz guided us in the definition of the insertion algorithm. Conversely, in type $B$, now that we defined an insertion algorithm, we can guess the corresponding matrix ansatz. Let $\lambda = (\lambda_n)_{n \geqslant 0} \in \mathbb{Q}(\alpha, \beta, \gamma, \delta, q, u, z)^{\mathbb{N}}$ and $g_\lambda : \{\bullet, \circ\}^* \to \mathbb{Q}(\alpha, \beta, \gamma, \delta, q, u, z)$ a function satisfying

$$\begin{cases} ug^\lambda(X \bullet \circ Y) & = qg^\lambda(X \circ \bullet Y) + \lambda_n(g^\lambda(X \bullet Y) + g^\lambda(X \circ Y)) \\ q^{n-1}(zq\alpha + u\delta)g^\lambda(\circ X) & = u^{n-1}(zq\beta + u\gamma)g^\lambda(\bullet X) \\ \beta g^\lambda(X \bullet) & = \delta g^\lambda(X \circ) + \lambda_n g^\lambda(X) \\ g^\lambda(\varepsilon) & = \lambda_0. \end{cases} \quad (4.1)$$

Our matrix ansatz generalise the matrix ansatz defined by Corteel, Josuat-Vergès and Williams [4] of the type $B$ case with $\gamma = \delta = 0$. As in section 2, these recurrences imply global recurrences but we do not know if the global recurrences imply the system (4.1). Nonetheless, the global recurrences we obtain are compatible, analogously to Theorem 9, with the insertion algorithm we defined.